# Privacy and Confidentiality in an e-Commerce World: Data Mining, Data Warehousing, Matching and Disclosure Limitation

**Stephen E. Fienberg**

*Abstract.* The growing expanse of e-commerce and the widespread availability of online databases raise many fears regarding loss of privacy and many statistical challenges. Even with encryption and other nominal forms of protection for individual databases, we still need to protect against the violation of privacy through linkages across multiple databases. These issues parallel those that have arisen and received some attention in the context of homeland security. Following the events of September 11, 2001, there has been heightened attention in the United States and elsewhere to the use of multiple government and private databases for the identification of possible perpetrators of future attacks, as well as an unprecedented expansion of federal government data mining activities, many involving databases containing personal information. We present an overview of some proposals that have surfaced for the search of multiple databases which supposedly do not compromise possible pledges of confidentiality to the individuals whose data are included. We also explore their link to the related literature on privacy-preserving data mining. In particular, we focus on the matching problem across databases and the concept of "selective revelation" and their confidentiality implications.

## 1. INTRODUCTION

Click on *Google* and search for "Feinberg contingency talb" and you will be asked if you meant "feinberg contingency table," and if you click on this again you will reach a mix of links to publications that refer to "Bishop, Fienberg and Holland" [3] or "Bishop, Feinberg and Holland," or other papers by the present author with his name spelled "Feinberg," "Fienberg" and many other ways! All thanks to the data mining tool of hidden Markov models and Google's page-rank methodology. This represents data mining at work in e-commerce, but in situations that do not violate my privacy or impinge on promises of confidentiality. Indeed, most authors in statistics are happy to have their name appear in a Google search whether it is spelled correctly or incorrectly. Data mining tools help enable searches as we engage in e-commerce, whether it is in a form like collaborative filtering or something more elaborate. When the data used by individual e-commerce vendors are linked to other databases, however, issues of privacy and confidentiality become front and center [29]. This has become of special concern in recent months as the U.S. government has attempted to secure individually identified information from Google and other companies engaged in e-commerce. (Katie Hafner and Matt Richtel, "Google Resists U.S. Subpoena of Search Data," *The New York Times,* January 20, 2006.)

*Stephen E. Fienberg is Maurice Falk University Professor of Statistics and Social Science, Carnegie Mellon University, Pittsburgh, Pennsylvania 15213, USA e-mail: fienberg@stat.cmu.edu.*

The website of the American Civil Liberties Union includes a "flash movie" of a telephone pizza order (www.aclu.org/pizza/) that triggers a series of data retrievals from some gigantic integrated database that includes medical records, travel information, magazine subscriptions, clothing purchases and seemingly instantaneously linked local area crime reports. It represents the public's worst fears regarding the invasion of privacy that has come from e-commerce and growth and spread of data warehousing. The website warns that "Government programs such as MATRIX and Carnivore are destroying our privacy. We live in a democratic society and government-controlled data systems are a dangerous step toward establishing a 24-hour surveillance society." What are these programs? Is the pizza movie myth or reality?

Here are some related stories in the news this past year:

- "Identity thieves posing as legitimate businesses were able to access profiles that include Social Security numbers, credit histories, criminal records and other sensitive material, ChoicePoint spokesman Chuck Jones said. [ChoicePoint] maintains personal profiles of nearly every U.S. consumer, which it sells to employers, landlords, marketing companies and about 35 U.S. government agencies. In California, the only state that requires companies to disclose security breaches, ChoicePoint sent warning letters to 30,000 to 35,000 consumers advising them to check their credit reports." Excerpted from "Thieves Steal Consumer Info Database," *CNN Money*, February 15, 2005.
- "Billions of records about virtually every adult in the country are maintained by an array of companies. Among the most familiar are the credit bureaus that have long tracked debts and payment histories. Less familiar, though, are data brokers such as ChoicePoint, which aggregate other personal information and operate with fewer restrictions. And, increasingly, banks and credit card companies maintain considerable data caches on their customers." Excerpted from "Firms Hit by ID Theft Find Way to Cash In on Victims," *Los Angeles Times*, August 22, 2005.
- "LexisNexis, a worldwide provider of legal and business data, announced yesterday that information about 32,000 consumers was fraudulently gathered in a series of incidents. The data include names, addresses and Social Security and driver's license numbers." Excerpted from Jonathan Krim and Robert O'Harrow Jr., "Data Under Siege," *The Washington Post*, Thursday, March 10, 2005 (www.washingtonpost.com/wp-dyn/articles/A19982-2005Mar9.html).

"Reed Elsevier, owner of the LexisNexis data bases, said Tuesday that Social Security numbers, driver's license information and the addresses of 310,000 people may have been stolen, 10 times more than it originally reported last month." Excerpted from "Security Breach at LexisNexis Now Appears Larger," by Heather Timmons, *The New York Times*, April 13, 2005.

- "...41 graduate students in a computer security course at Johns Hopkins University...became mini-data-brokers themselves over the last semester.... Working with a strict requirement to use only legal, public sources of information, groups of three to four students set out to vacuum up not just tidbits on citizens of Baltimore, but whole databases: death records, property tax information, campaign donations, occupational license registries. They then cleaned and linked the databases they had collected, making it possible to enter a single name and generate multiple layers of information on individuals. Each group could spend no more than $50. ...Several groups managed to gather well over a million records, with hundreds of thousands of individuals represented in each database." Excerpted from "Personal Data for the Taking," by Tom Zeller Jr., *The New York Times*, May 18, 2005.
- "In one of the largest breaches of data security to date, CitiFinancial, the consumer finance subsidiary of Citigroup, announced yesterday that a box of computer tapes containing information on 3.9 million customers was lost by United Parcel Service last month, while in transit to a credit reporting agency." Excerpted from "Personal Data for 3.9 Million Lost in Transit," by Tom Zeller Jr., *The New York Times*, June 7, 2005.

Data warehousing companies such as Acxiom, ChoicePoint and LexisNexis use their data to perform background checks on prospective applicants to employers, insurers and credit providers. They also sell their data to state and federal governments. Figure 1 shows the array of data available from ChoicePoint and the types of clients who access, as presented by *The Washington Post* (www.washingtonpost.com/wp-srv/business/daily/graphics/choicepoint_012005.html). If you go to the ChoicePoint website (www.choicepoint.com) and read the privacy policy you are told about "How we protect you," but if you want to check the accuracy of information on yourself that ChoicePoint sells to others you need to provide your Social Security number! This means that if ChoicePoint did not have your Social Security number before, it would now, and they make no promise about how it will (or will not) be used or shared in the future.

In part as a consequence of the data security breaches of the sort described above, some form of data breach legislation has been introduced in at least 35 states and signed into law in at least 22, according to data compiled by the National Conference of State Legislatures. (Tom Zeller Jr. "Link by link; Waking up to recurring ID nightmares," *The New York Times*, January 9, 2006.)

In the next section we briefly describe a related set of government data mining and data warehousing activities that came into the public eye following the terrorist attacks of September 11, 2001. The link with the more public e-commerce activities was MATRIX, referred to by the ACLU webpage but which has since been "publicly" abandoned. In Section 3 we give an overview of record linkage and its use for merging large data files from diverse sources as well as its implications for the splitting of databases for privacy protection. Section 4 reviews some proposals that have surfaced for the search of multiple databases without compromising possible pledges of confidentiality to the individuals whose data are included and their link to the related literature on privacy-preserving data mining. In particular, we focus on the concept of *selective revelation* and its confidentiality implications. We relate these ideas to the recent statistical literature on disclosure limitation for confidential databases and explain the problems with the privacy claims. We conclude with some observations regarding privacy protection and e-commerce.

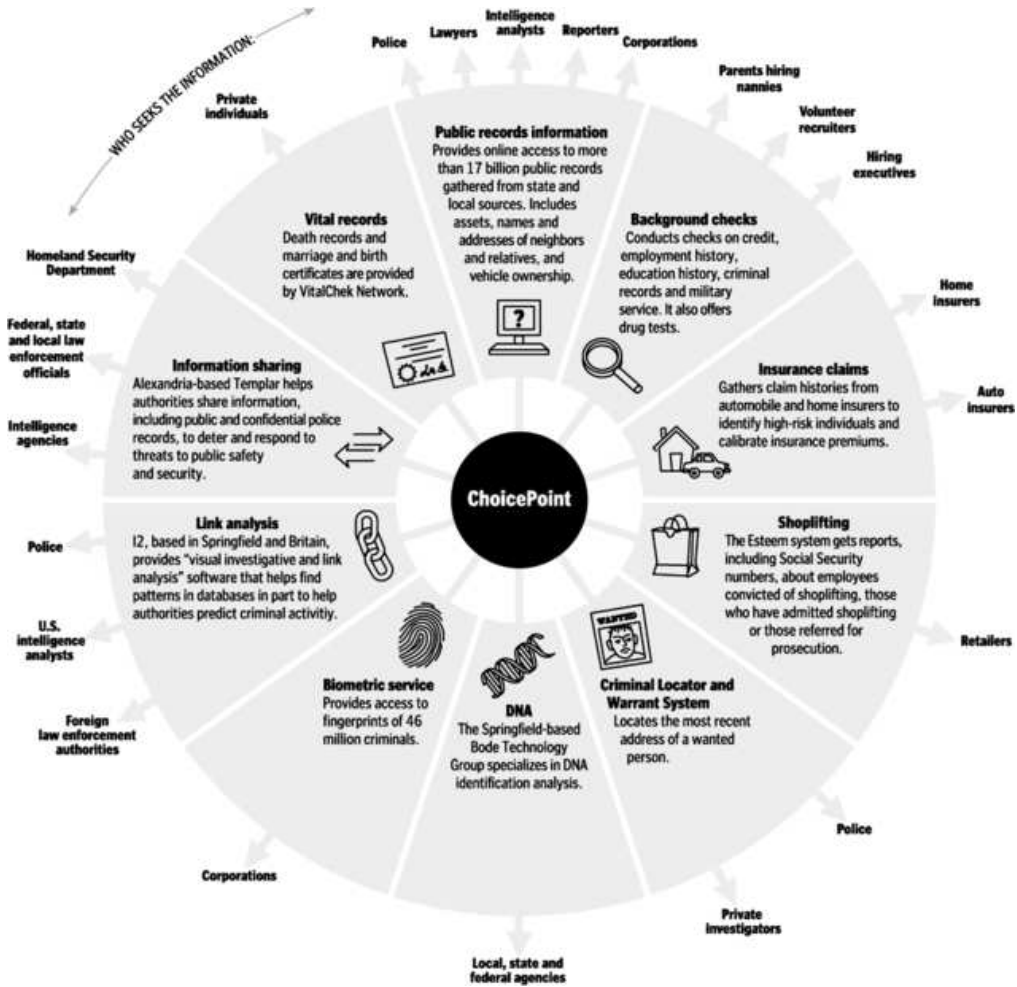

Fig. 1.   *ChoicePoint data sources and clients. Source: The Washington Post, January 20, 2005.*

## 2. HOMELAND SECURITY AND THE SEARCH FOR TERRORISTS

A recently issued report from the U.S. General Accounting Office [37] notes that at least 52 agencies are using or planning to use data mining, "factual data analysis," or "predictive analytics," in some 199 different efforts. Of these, at least 29 projects involve analyzing intelligence and detecting terrorist activities, or detecting criminal activities or patterns. Notable among the nonresponders to the GAO inquiry were agencies like the Central Intelligence Agency and the National Security Agency (NSA).

Perhaps the most visible of these efforts was the *Total Information Awareness* (TIA) program initiated by the Defense Advanced Research Program (DARPA) in DARPA's Information Awareness Office (IAO), which was established in January 2002, in the aftermath of the September 11 terrorist attacks. The TIA research and development program was aimed at integrating information technologies into a prototype to provide tools to better detect, classify and identify potential foreign terrorists. When it came under public scrutiny in 2003, TIA morphed into the *Terrorist Information Program* (still TIA) with essentially the same objectives, although it too did not move forward into implementation. TIA served as the model, however, for the *Multi-state Anti-terrorism Information Exchange* system (MATRIX) that was in use in seven states for a period of time during 2004 and 2005, and was intended to provide "the capability to store, analyze, and exchange sensitive terrorism-related information in MATRIX data bases among agencies, within a state, among states, and between state and federal agencies."

According to a recent report from the Congressional Research Service [30] [footnotes omitted]:

> The MATRIX project was initially developed in the days following the September 11, 2001, terrorist attacks by Seisint, a Florida-based information products company, in an effort to facilitate collaborative information sharing and factual data analysis. At the outset of the project, MATRIX included a component Seisint called the High Terrorist Factor (HTF), which was designed to identify individuals with high HTF scores, or so-called terrorism quotients, based on an analysis of demographic and behavioral data. Although the HTF scoring system appeared to attract the interest of officials, this feature was reportedly dropped from MATRIX because it relied on intelligence data not normally available to the law enforcement community and because of concerns about privacy abuses.
> ...The analytical core of the MATRIX pilot project is an application called Factual Analysis Criminal Threat Solution (FACTS), described as a "technological, investigative tool allowing query-based searches of available state and public records in the data reference repository." The FACTS application allows an authorized user to search "dynamically" combined records from disparate datasets based on partial information, and will "assemble" the results. The data reference repository used with FACTS represents the amalgamation of over 3.9 billion public records collected from thousands of sources. The data contained in FACTS include FAA pilot license and aircraft ownership records, property ownership records, information on vessels registered with the Coast Guard, state sexual offender lists, federal terrorist watch lists, corporation filings, Uniform Commercial Code filings, bankruptcy filings, state-issued professional license records, criminal history information, department of corrections information and photo images, driver's license information and photo images, motor vehicle registration information, and information from commercial sources that "are generally available to the public or legally permissible under federal law."
> ...To help address the privacy concerns associated with a centralized data repository, some officials have suggested switching to a distributed approach whereby each state would maintain possession of its data and control access according to its individual laws.

The data reference repository is said to exclude data from the following sources:

- telemarketing call lists,
- direct mail mailing lists,

- airline reservations or travel records,
- frequent flyer/hotel stay program membership information or activity,
- magazine subscription records,
- information about purchases made at retailers or over the Internet,
- telephone calling logs or records,
- credit or debit card numbers,
- mortgage or car payment information,
- bank account numbers or balance information,
- records of birth certificates, marriage licenses and divorce decrees, and
- utility bill payment information.

Nonetheless, MATRIX and its data records sound suspiciously like the ACLU Pizza Movie scenario! And the links to the news stories excerpted in Section 1 are more direct than one might imagine. In 2004, LexisNexis acquired Seisint and the security breaches were in the new Seisint subsidiary, the very same one that provides the data for MATRIX!

MATRIX was officially abandoned as a multistate activity in April 2005 although individual states were allowed to continue with their parts of the program. This does not mean the demise of the TIA effort, however, as there are other federal initiatives built on a similar model:

- Analysis, Dissemination, Visualization, Insight and Semantic Enhancement (ADVISE), which is a research and development program within the Department of Homeland Security (DHS), part of its three-year-old "Threat and Vulnerability, Testing and Assessment" portfolio (Mark Clayton, "US Plans Massive Data Sweep," *The Christian Science Monitor*, February 9, 2006. www.csmonitor.com/2006/0209/p01s02-uspo.html).
- The Information Awareness Prototype System (IAPS), the core architecture that tied together numerous information extraction, analysis and dissemination tools developed under TIA, including the privacy-protection technologies, was moved to the Advanced Research and Development Activity (ARDA), housed at NSA headquarters in Fort Meade, Md (Shane Harris, "TIA Lives On," *National Journal*, Thursday, Feb. 23, 2006).

In TIA, MATRIX, ADVISE and IAPS, the data miner can issue queries to the multiple linked databases and receive responses that combine data on individuals across the databases. The goal is the identification of terrorists or criminals in a way that would not be possible from the individual databases. We distinguish between two aspects of this goal: (1) identification of known terrorists which is a form of retro- or postdiction, and (2) identification of potential future terrorists and profiling, which involves prediction. Prediction cannot be separated from uncertainty; postdiction might conceivably be. Most of the public outcry regarding TIA and MATRIX has focused on concerns regarding what has been described as "dataveillance" [4] and terrorist profiling, that is, concerns both about the use of data for purposes other than those for which they were collected without the consent of the individual, and about the quality and accuracy of the mined data and the likelihood that they may help falsely identify individuals as terrorists.

In the next two sections, we explore some issues related to the creation and the use of "linked" databases for the privacy of the individuals whose confidential information is contained in them.

## 3. MATCHING AND RECORD LINKAGE METHODS

More than 100 vendors offer record matching systems, some of which sell for thousands of dollars, but most of the underlying methodology for such systems is proprietary and few details are publicly available. Matches can occur at random. For example, consider a pair of files, $A$ and $B$, containing $n$ records on the same individuals. Then the probability of correctly matching exactly $r$ individuals by picking a random permutation for file $B$ and linking to file $A$ is

$$\frac{\sum_{\nu=0}^{n-r}((-1)^n - r)/\nu!}{r!}. \tag{1}$$

Domingo-Ferrer and Torra [8] derive this baseline and illustrate it numerically in an example with $n = 90$, where the expected number of correct matches is $O(10^{24})$. Working with actual data in the matching process can change this situation drastically.

Bilenko et al. [2] provide an overview of the published literature on the topic noting that most methods rely on the existence of unique identifiers or use some variation of the algorithm presented in Fellegi and Sunter [14]. Fellegi and Sunter's approach is built on several key components for identifying matching pairs of records across two files:

- Represent every pair of records using a vector of features (variables) that describe similarity between individual record fields. Features can be Boolean (e.g., last-namematches), discrete (e.g., first-$n$-characters-of-name-agree) or continuous (e.g., string-edit-distance-between-first-names).
- Place feature vectors for record pairs into three classes: matches ($M$), nonmatches ($U$) and possible matches. These correspond to "equivalent," "nonequivalent" and possibly equivalent (e.g., requiring human review) record pairs, respectively.
- Perform record-pair classification by calculating the ratio $(P(\gamma \mid M))/(P(\gamma \mid U))$ for each candidate record pair, where $\gamma$ is a feature vector for the pair and $P(\gamma \mid M)$ and $P(\gamma \mid U)$ are the probabilities of observing that feature vector for a matched and nonmatched pair, respectively. Two thresholds based on desired error levels—$T_\mu$ and $T_\lambda$—optimally separate the ratio values for equivalent, possibly equivalent and nonequivalent record pairs.
- When no training data in the form of duplicate and nonduplicate record pairs is available, matching can be unsupervised, where conditional probabilities for feature values are estimated using observed frequencies.
- Because most record pairs are clearly nonmatches, we need not consider them for matching. The way to manage this is to "block" the databases, for example, based on geography or some other variable in both databases, so that only records in comparable blocks are compared. Such a strategy significantly improves efficiency.

The first four components lay the groundwork for accuracy of record-pair matching using statistical techniques such as logistic regression, the EM algorithm and Bayes networks (e.g., see [22, 25, 38]). Accuracy is well known to be high when there is a 1–1 match between records in the two systems and deteriorates as the overlap between the files decreases as well as with the extent of measurement error in the feature values. While the use of human review of possible matches has been an integral part of many statistical applications, it may well be infeasible for large-scale data warehousing. The fifth component provides for efficiently processing large databases, but to the extent that blocking is approximate and possibly inaccurate its use decreases the accuracy of record-pair matching.

There are three potential lessons associated with this literature on matching and the methods it has produced:

1. If we are trying to protect against an intruder who would like to merge the data in a confidential database with an external database in his/her possession, then we need to assure ourselves and the intruder that the accuracy of matching is low and that individuals cannot be identified with high probability. We need to keep in mind that an intruder will have easy access to a host of identifiable public record systems. For example, as of September 7, 2005, *SearchSystems.net* (www.searchsystems.net/) listed 34,035 free searchable public record databases on its website!
2. One strategy for protecting a database against attack from an intruder is to split it into parts, perhaps overlapping, to decrease the likelihood of accurate matches. The parts should be immune from attack (with high probability) but of value for analytical purposes. For categorical data this might correspond to reporting lower-dimensional margins from a high-dimensional contingency table; see [5, 6] and [18]. For continuous data we might need to apply disclosure protection methods to the split components; for example, see [9] and [15] for overviews. It is the uncertainty associated with efforts to concatenate the separate pieces that provides the confidentiality protection in both instances. The higher the uncertainty the better the protection.
3. Unless ChoicePoint and other data warehousers are adding data into their files using unique identifiers such as Social Security numbers (and even Social Security numbers are not really unique!), or with highly accurate addresses and/or geography, some reasonable fraction of the data in their files will be the result of inaccurate and faulty matches. Data quality for data warehouses is an issue we all need to worry about; see [39].

## 4. ENCRYPTION, MULTIPARTY COMPUTATION AND PRIVACY-PRESERVING DATA MINING

If you search the WWW for "e-commerce" and "data privacy protection" you will find extensive discussion about firewalls, intrusion prevention (IPS) and intrusion detection (IDS) systems, and secure socket layer (SSL) encryption technology. Indeed, these technological tools are important for secure data transmission, statistical production and offline data storage; see [7]. But encryption cannot protect the privacy of individuals whose data are available in online databases!

Among the methods advocated to carry out such data mining exercises are those that are described as privacy-preserving data mining (PPDM). PPDM typically refers to data mining computations performed on the combined data sets of multiple parties without revealing each party's data to the other parties. The data consist of possibly overlapping sets of variables contained in the separate databases of the parties and overlapping sets of individuals. When the parties have data for the same variables but different individuals the data are said to be horizontally partitioned, whereas when the individuals are the same but the variables are different the data are said to be vertically partitioned. Here we are concerned with the more complex case involving both overlapping variables *and* overlapping sets of individuals. PPDM research comes in two varieties. In the first, sometimes referred to as the construction of "privacy-preserving statistical databases," the data are altered prior to delivery for data mining, for example, through the addition of random noise or some other form of perturbation. While these approaches share much in common with the methods in the literature on statistical disclosure limitation, they are of little use when it comes to the identification of terrorists. In the second variety, the problem is solved using what is known as "multiparty secure computation," where no party knows anything except its own input and the results. The literature typically presumes that data are included without error and thus could be matched perfectly if only there were no privacy concerns. The methods also focus largely on situations where the results are of some computation, such as a dot product or the description of an association rule. See the related discussion in [19].

A major problem with the PPDM literature involving multiparty computation is that the so-called proofs of security are designed to protect not the individuals in the database but rather the database owners, as in the case of two companies sharing information but not wanting to reveal information about their customers to one another beyond that contained in the shared computation. Once the results of the data mining consist of linked extracts of the data themselves, however, the real question is whether one of the parties can use the extra information to infer something about the individuals in the other party's data that would otherwise not be available.

Secure computation is a technique for carrying out computations across multiple databases without revealing any information about data elements found only in one database. The technique consists of a protocol for exchanging messages. We assume the parties to be *semihonest*: that is, they correctly follow the protocol specification, yet attempt to learn additional information by analyzing the messages that are passed. For example, Agrawal, Evfimievski and Srikant [1] illustrate the secure computation notion via an approach to the matching problem for parties $A$ and $B$. They introduce a pair of encryption functions $E$ (known only to $A$) and $E'$ (known only to $B$) such that for all $x$, $E(E'(x)) = E'(E(x))$. $A$'s database consists of a list $\mathbf{A}$ and $B$'s consists of a list $\mathbf{B}$. $A$ sends $B$ the message $E(\mathbf{A})$; $B$ computes $E'(E(\mathbf{A}))$ and then sends to $A$ the two messages $E'(E(\mathbf{A}))$ and $E'(\mathbf{B})$. $A$ then applies $E$ to $E'(\mathbf{B})$, yielding $E'(E(A))$ and $E'(E(B))$. $A$ computes $E'(E(\mathbf{A})) \cap E'(E(\mathbf{B}))$. Since $A$ knows the order of items in $\mathbf{A}$, $A$ also knows the order of items in $E'E(A))$ and can quickly determine $\mathbf{A} \cap \mathbf{B}$. The main problems with this approach are (1) it is asymmetric, that is, $B$ must trust $A$ to send $\mathbf{A} \cap \mathbf{B}$ back, and (2) it presumes semihonest behavior.

Li, Tygar and Hellerstein in [26] describe a variety of scenarios in which the Agrawal et al. protocol can easily be exploited by one party to obtain a great deal of information about the other's database, and they explain drawbacks of some other secure computation methods including the use of one-way hash-based schemes. As Dwork and Nissim [13] note: "There is also a very large literature in secure multi-party computation. In secure multi-party computation, functionality is paramount, and privacy is only preserved to the extent that the function outcome itself does not reveal information about the individual inputs. In privacy-preserving statistical data bases, privacy is paramount." The problem with privacy-preserving datamining methods for terrorist detection is that they seek the protection of the latter while revealing individual records using the functionality of the former. For more details on some of these and other issues, see [23].

The U.S. Congress and various private foundations have taken up the issue of privacy protection from government data mining activities especially in the

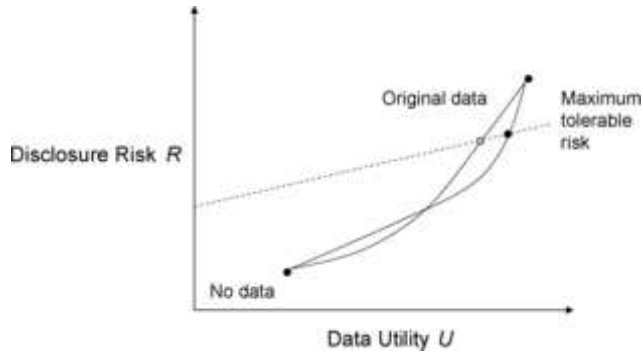

Fig. 2. *R–U confidentiality maps for two different disclosure limitation methods with varying parameter settings. Adapted from* [12].

post-9/11 world. For example, in its recent report, the U.S. Department of Defense Technology and Privacy Advisory Committee (TAPAC) [36] has stressed the existence of a broad array of government data mining programs, and disjointed, inconsistent and outdated laws and regulations protecting privacy. TAPAC recommended broad new actions to protect privacy, both within the Department of Defense and across agencies of the federal government.

The long-standing concern regarding surveillance of U.S. citizens and others by government agencies has been heightened during the war on terror (e.g., see [24]) and especially most recently with the controversy over unauthorized domestic spying. (David Johnston and Neil A. Lewis, "Domestic Surveillance: The White House; Defending Spy Program, Administration Cites Law," *The New York Times*, December 23, 2005.)

## 5. SELECTIVE REVELATION, THE RISK-UTILITY TRADE-OFF AND DISCLOSURE LIMITATION ASSESSMENT

To get around the privacy problems associated with the development of the TIA and MATRIX systems Tygar [34, 35] and others have advocated the use of what has come to be called "selective revelation," involving something like the risk-utility trade-off in statistical disclosure limitation. Sweeney [33] used the term to describe an approach to disclosure limitation that allows data to be shared for surveillance purposes "with a sliding scale of identifiability, where the level of anonymity matches scientific and evidentiary need." This corresponds to a monotonically increasing threshold for maximum tolerable risk in the R–U confidentiality map framework described in [10, 11, 12], as depicted in Figure 2.

Figure 3 depicts the basic selective revelation scheme as described in a committee report on TIA privacy methodology [21].

The TIA privacy report [21] suggests that

> Selective revelation works by putting a security barrier between the private data and the analyst, and controlling what information can flow across that barrier to the analyst. The analyst injects a query that uses the private data to determine a result, which is a high-level sanitized description of the query result. That result must not leak any private information to the analyst. Selective revelation must accommodate multiple data sources, all of which lie behind the (conceptual) security barrier. Private information is not made available directly to the analyst, but only through the security barrier.

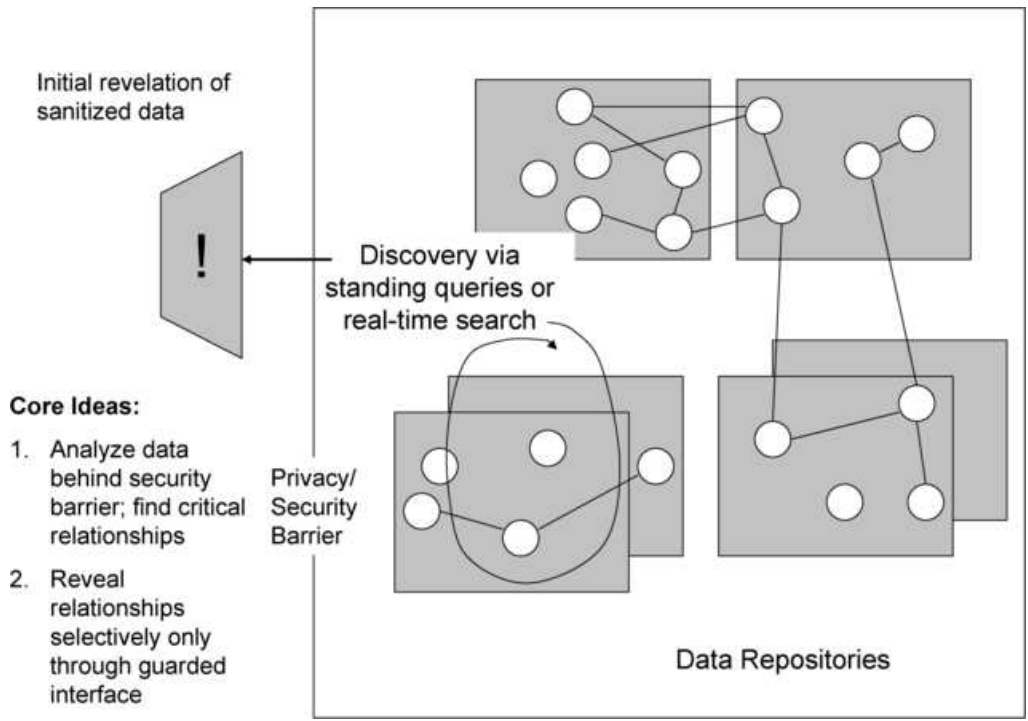

FIG. 3. *Idealized selective revelation architecture. Adapted from Slide* 11, [21].

One effort to implement this scheme was dubbed *privacy appliances* by Lunt [27] and it was intended to be a stand-alone device that would sit between the analyst and the private data source so that private data stays in authorized hands. These privacy controls would also be independently operated to keep them isolated from the government. According to Lunt [27] the device would provide:

- *Inference control* to prevent unauthorized individuals from completing queries that would allow identification of ordinary citizens.
- *Access control* to return sensitive identifying data only to authorized users.
- *Immutable audit trail* for accountability.

Implicit in the TIA Report and in the Lunt approach was the notion that linkages across databases behind the security barrier would utilize identifiable records and thus some form of multiparty computation method involving encryption techniques.

The real questions of interest in "inference control" are: (1) What disclosure limitation methods should be used? (2) To which databases should they be applied? and (3) How can the "inference control" approaches be combined with the multiparty computation methods? Here is what we know in the way of answers:

1. Both Sweeney [33] and Lunt et al. [28] refer to Sweeney's version of microaggregation, known as $k$-anonymity, but with few details on how it could be used in this context. This methodology combines observations in groups of size $k$ and reports either the sum or the average of the group for each unit. The groups may be identified by clustering or some other statistical approach. Left unsaid is what kinds of analyses users might perform with such aggregated data. Further, neither $k$-anonymity nor any other confidentiality tool does anything to cope with the implications of the release of exactly linked files requested by "authorized users."

2. Much of the statistical and operations research literature on confidentiality fails to address the risk-utility trade-off, largely by focusing primarily only on privacy, or on technical implementations without understanding how users wish to analyze a database; for example, see [20].
3. A clear lesson from the statistical disclosure limitation literature is that privacy protection in the form of "safe releases" from separate databases does not guarantee privacy protection for a merged database. A figure in [28] demonstrates recognition of this fact by showing privacy appliances applied for the individual databases and then, again, independently for the combined data.
4. To date there have been a limited number of crosswalks between the statistical disclosure limitation literatures on multiparty computation and risk-utility trade-off choices for disclosure limitation. Zhong, Yang and Wright [40] provide a starting point for discussions on $k$-anonymity. There are clearly a number of alternatives to $k$-anonymity, and ones which yield "anonymized" databases of far greater statistical utility!
5. The hype associated with the TIA approach to protection has abated, largely because TIA no longer exists as an official program. But similar programs continue to appear in different places in the federal government and no one associated with any of them has publicly addressed the privacy concerns raised here regarding the TIA approach.

When the U.S. Congress stopped the funding for DARPA's TIA program in 2003, Lunt's research and development effort at PARC Research Center was an attendant casualty. Thus to date there have been no publicly available prototypes of the privacy appliance, nor are there likely to be in the near future. The claims of privacy protection and selective revelation continue with MATRIX and other data warehouse systems, but without an attendant research program, and the federal government continues to plan for the use of data mining techniques in other federal initiatives such as the Computer Assisted Passenger Profiling System II (CAPPS II). Similar issues arise in the use of government, medical and private transactional data in bio-terrorism surveillance; for example, see [17] and [32].

## 6. CONCLUSIONS

Data privacy protection is a major issue for e-commerce. While solutions like SSL encryption may help companies with protection for confidential data transmission, the privacy pitfalls of marketing data as part of e-commerce are many. In this paper, we have focused on large-scale data warehousing in part because the repeated announcements of security breaches in systems operated by the major vendors such as Acxiom, ChoicePoint and LexusNexus have filled our morning newspapers during the past several years. The public and civil rights groups have argued that this is just the tip of the privacy-violation iceberg and they have called for government intervention and legal restrictions on both public and private organizations with respect to data warehousing and data mining. The lessons from such privacy breaches extend easily to virtually all electronically accessible databases. Companies need to take data security seriously and implement "best practices," and they need to rethink their policies on "data access" by others.

The giant data warehouses described in this paper have been assembled through the aggregation of information from many separate databases and transactional data systems. They depend heavily on matching and record-linkage methods that intrinsically are statistical in nature, and whose accuracy deteriorates rapidly in the presence of serious measurement error. Data mining tools cannot make up for bad data and poor matches, and someone beyond "wronged consumers" will soon begin to pay attention.

Should you worry about these data warehouses? With very high probability they contain data on you and your household, but you will never quite know what data or how accurate the information is. And soon the data may be matched into government-sponsored terrorist search systems such as the one being set up by the Transportation Security Administration (TSA) to match passenger lists into a consolidated watch list of suspected terrorists. On September 19, 2005, the Secure Flight Working Group to the Transportation Security Administration (TSA) submitted a report questioning TSA's secrecy regarding what data it plans to use and how [31]:

> The TSA is under a Congressional mandate to match domestic airline passenger lists against the consolidated terrorist watch list. TSA has failed to specify with consistency whether watch list matching is the only goal of Secure Flight at this stage...
> Will Secure Flight be linked to other TSA applications?...
> How will commercial data sources be used? One of the most controversial elements of Secure Flight has been the possible uses of commercial data. TSA has never clearly defined two threshold issues: what it means by "commercial data"; and how it might use commercial data sources in the implementation of Secure Flight. TSA has never clearly distinguished among various possible uses of commercial data, which all have different implications.

The story continues, however, since a few months later it was revealed that TSA had purchased a database from ChoicePoint to be matched against the watch list. ("TSA Chief Suspends Traveler Registry Plans," Associated Press, February 9, 2006.)

Finally, we need new computational and statistical technologies to protect linked multiple databases from privacy protection in the face of commercial and government queries. Slogans like "selective revelation" are not enough without technical backup. This might be provided by the serious integration of research ideas emanating from the statistical disclosure and cryptography communities. The technologies that result from such collaborative research must be part of the public domain, because only then can we evaluate their adequacy.

## ACKNOWLEDGMENTS

The research reported here was supported in part by NSF Grants EIA-98-76619 and IIS-01-31884 to the National Institute of Statistical Sciences and by Army Contract DAAD19-02-1-3-0389 to CyLab at Carnegie Mellon University. This paper is based in part on an earlier and much shorter paper focusing on homeland security issues; see [16]. I have benefited from conversations with Chris Clifton, Cynthia Dwork, Alan Karr and Latanya Sweeney about the material described here but they bear no responsibility for how I have represented their input. I thank the referees for comments that improved the content of the paper.

## REFERENCES

[1] AGRAWAL, R., EVFIMIEVSKI, A. and SRIKANT, R. (2003). Information sharing across private databases. In *Proc. 2003 ACM SIGMOD International Conference on Management of Data* 86–97. ACM Press, New York.

[2] BILENKO, M., MOONEY, R., COHEN, W. W., RAVIKUMAR, P. and FIENBERG, S. E. (2003). Adaptive name matching in information integration. *IEEE Intelligent Systems* **18**(5) 16–23.
[3] BISHOP, Y. M. M., FIENBERG, S. E. and HOLLAND, P. W. (1975). *Discrete Multivariate Analysis*: *Theory and Practice*. MIT Press, Cambridge, MA. MR0381130
[4] CLARKE, R. (1988). Information technology and dataveillance. *Comm. ACM* **31** 498–512.
[5] DOBRA, A. and FIENBERG, S. E. (2001). Bounds for cell entries in contingency tables induced by fixed marginal totals. *Statist. J. United Nations ECE* **18** 363–371.
[6] DOBRA, A. and FIENBERG, S. E. (2003). Bounding entries in multi-way contingency tables given a set of marginal totals. In *Foundations of Statistical Inference* (Y. Haitovsky, H. R. Lerche and Y. Ritov, eds.) 3–16. Physica, Heidelberg. MR2017809
[7] DOMINGO-FERRER, J., MATEO-SANZ, J. M. and SÁNCHEZ DEL CASTILLO, R. X. (2000). Cryptographic techniques in statistical data protection. In *Proc. Joint UN/ECE-Eurostat Work Session on Statistical Data Confidentiality* 159–166. Office for Official Publications of the European Communities, Luxembourg.
[8] DOMINGO-FERRER, J. and TORRA, V. (2003). Disclosure risk assessment in statistical microdata protection via advanced record linkage. *Stat. Comput.* **13** 343–354. MR2005437
[9] DUNCAN, G. T. (2001). Confidentiality and statistical disclosure limitation. *International Encyclopedia of the Social and Behavioral Sciences* 2521–2525. North-Holland, Amsterdam.
[10] DUNCAN, G. T., FIENBERG, S. E., KRISHNAN, R., PADMAN, R. and ROEHRIG, S. F. (2001). Disclosure limitation methods and information loss for tabular data. In *Confidentiality, Disclosure and Data Access*: *Theory and Practical Applications for Statistical Agencies* (P. Doyle, J. Lane, J. Theeuwes and L. Zayatz, eds.) 135–166. North-Holland, Amsterdam.
[11] DUNCAN, G. T., KELLER-MCNULTY, S. A. and STOKES, S. L. (2004). Database security and confidentiality: Examining disclosure risk vs. data utility through the R–U confidentiality map. Technical Report 142, National Institute of Statistical Sciences.
[12] DUNCAN, G. T. and STOKES, S. L. (2004). Disclosure risk vs. data utility: The R–U confidentiality map as applied to topcoding. *Chance* **17**(3) 16–20. MR2061932
[13] DWORK, C. and NISSIM, K. (2004). Privacy-preserving data mining on vertically partitioned databases. In *Proc. CRYPTO 2004*, *24th International Conference on Cryptology* 528–544. Univ. California, Santa Barbara.
[14] FELLEGI, I. P. and SUNTER, A. B. (1969). A theory for record linkage. *J. Amer. Statist. Assoc.* **64** 1183–1210.
[15] FIENBERG, S. E. (2005). Confidentiality and disclosure limitation. *Encyclopedia of Social Measurement* 463–469. North-Holland, Amsterdam.
[16] FIENBERG, S. E. (2005). Homeland insecurity: Datamining, terrorism detection, and confidentiality. *Bull. Internat. Stat. Inst.*, 55th Session. Sydney.
[17] FIENBERG, S. E. and SHMUELI, G. (2005). Statistical issues and challenges associated with rapid detection of bio-terrorist attacks. *Stat. Med.* **24** 513–529. MR2134521
[18] FIENBERG, S. E. and SLAVKOVIC, A. B. (2004). Making the release of confidential data from multi-way tables count. *Chance* **17**(3) 5–10. MR2061930
[19] FIENBERG, S. E. and SLAVKOVIC, A. B. (2005). Preserving the confidentiality of categorical statistical data bases when releasing information for association rules. *Data Mining and Knowledge Discovery* **11** 155–180.
[20] GOPAL, R., GARFINKEL, R. and GOES, P. (2002). Confidentiality via camouflage: The CVC approach to disclosure limitation when answering queries to databases. *Oper. Res.* **50** 501–516. MR1910286
[21] INFORMATION SCIENCE AND TECHNOLOGY STUDY GROUP ON SECURITY AND PRIVACY (chair: J. D. Tygar) (2002). Security With Privacy. Briefing.
[22] JARO, M. A. (1995). Probabilistic linkage of large public health data files. *Stat. Med.* **14** 491–498.
[23] KARR, A. F., LIN, X., SANIL, A. P. and REITER, J. P. (2006). Secure statistical analysis of distributed databases. In *Statistical Methods in Counterterrorism* (A. Wilson, G. Wilson and D. H. Olwell, eds.). Springer, New York.
[24] KREIMER, S. F. (2004). Watching the watchers: Surveillance, transparency, and political freedom in the war on terror. *J. Constitutional Law* **7** 133–181.
[25] LARSEN, M. D. and RUBIN, D. B. (2001). Iterative automated record linkage using mixture models. *J. Amer. Statist. Assoc.* **96** 32–41.

[26] LI, Y., TYGAR, J. D. and HELLERSTEIN, J. M. (2005). Private matching. In *Computer Security in the 21st Century* (D. T. Lee, S. P. Shieh and J. D. Tygar, eds.) 25–50. Springer, New York.

[27] LUNT, T. (2003). Protecting privacy in terrorist tracking applications. Presentation to the Department of Defense Technology and Privacy Advisory Committee, September 29, 2003.

[28] LUNT, T., STADDON, J., BALFANZ, D., DURFEE, G., URIBE, T. et al. (2005). Protecting privacy in terrorist tracking applications. Powerpoint presentation. Available at research.microsoft.com/projects/SWSecInstitute/five-minute/Balfanz5.ppt.

[29] MURALIDHAR, K., SARATHY, R. and PARSA, R. (2001). An improved security requirement for data perturbation with implications for e-commerce. *Decision Sci.* **32** 683–698.

[30] RELYEA, H. C. and SEIFERT, J. W. (2005). Information Sharing for Homeland Security: A Brief Overview. Congressional Research Service, The Library of Congress (Updated January 10, 2005). Available at `www.fas.org/sgp/crs/RL32597.pdf`.

[31] SECURE FLIGHT WORKING GROUP (2005). Report of the secure flight working group. Presented to the Transportation Security Administration, September 19, 2005. Available at www.epic.org/privacy/airtravel/sfwg_report_091905.pdf.

[32] SWEENEY, L. (2005). Privacy-preserving bio-terrorism surveillance. Presentation at AAAI Spring Symposium, AI Technologies for Homeland Security, Stanford Univ.

[33] SWEENEY, L. (2005). Privacy-preserving surveillance using selective revelation. LIDAP Working Paper 15, School Computer Science, Carnegie Mellon Univ.

[34] TYGAR, J. D. (2003). Privacy architectures. Presentation at Microsoft Research, June 18, 2003. Available at research.microsoft.com/projects/SWSecInstitute/slides/Tygar. pdf.

[35] TYGAR, J. D. (2003). Privacy in sensor webs and distributed information systems. In *Software Security Theories and Systems* (M. Okada, B. Pierce, A. Scedrov, H. Tokuda and A. Yonezawa, eds.) 84–95. Springer, New York.

[36] U.S. DEPARTMENT OF DEFENSE TECHNOLOGY AND PRIVACY ADVISORY COMMITTEE (TAPAC) (2004). *Safeguarding Privacy in the Fight Against Terrorism.* Department of Defense, Washington.

[37] U.S. GENERAL ACCOUNTING OFFICE (2004). *Data Mining*: *Federal Efforts Cover a Wide Range of Uses.* GAO-04-548, Report to the Ranking Minority Member, Subcommittee on Financial Management, the Budget and International Security, Committee on Governmental Affairs, U.S. Senate, Washington.

[38] WINKLER, W. E. (2002). Methods for record linkage and Bayesian networks. *Proc. Section Survey Research Methods* 3743–3748. Amer. Statist. Assoc., Alexandria, VA.

[39] WINKLER, W. E. (2005). Data quality in data warehouses. *Encyclopedia of Data Warehousing and Data Mining* **1**. Idea Group, Hershey, PA.

[40] ZHONG, S., YANG, Z. and WRIGHT, R. N. (2005). Privacy-enhancing $k$-anonymization of customer data. In *Proc. 24th ACM SIGMOD International Conference on Management of Data/Principles of Database Systems* (*PODS 2005*). ACM Press, New York.